\documentclass[letterpaper, 10 pt, conference]{ieeeconf}  

\IEEEoverridecommandlockouts                              

\usepackage[ansinew]{inputenc}
\usepackage{float}
\usepackage{graphicx}
\usepackage{fancyhdr}
\usepackage{amsmath,amssymb}
\usepackage{color}
\usepackage{verbatim}
\usepackage{scrtime}
\usepackage{float}
\usepackage{subfigure}
\usepackage{url}
\usepackage{multirow}
\usepackage[usenames,dvipsnames]{xcolor}

\usepackage{hyperref}
\hypersetup{
     colorlinks  = true,
     citecolor   = blue
}
\hypersetup{linkcolor=blue}

\usepackage{enumerate}

\usepackage{algpseudocode}
\usepackage{algorithm}
\algtext*{EndIf}
\algnewcommand\algorithmicinput{\textbf{Input:}}
\algnewcommand\Input{\item[\algorithmicinput]}
\algnewcommand\algorithmicoutput{\textbf{Ouput:}}
\algnewcommand\Output{\item[\algorithmicoutput]}

\title{\LARGE \bf
A Riemannian approach to low-rank algebraic Riccati equations*
}

\author{Bamdev Mishra$^{1}$ and Bart Vandereycken$^{2}$
\thanks{*This paper presents research results of the Belgian Network DYSCO (Dynamical Systems, Control, and Optimization), funded by the Interuniversity Attraction Poles Programme, initiated by the Belgian State, Science Policy Office. The scientific responsibility rests with its authors. This work was also partly supported by the Belgian FRFC (Fonds de la Recherche Fondamentale Collective). Bamdev Mishra is a research fellow of the Belgian National Fund for Scientific Research (FNRS).}
\thanks{$^{1}$Department of Electrical Engineering and Computer Science,
        University of Li\`ege, Li\`ege, 4000 Belgium
        {\tt\small B.Mishra@ulg.ac.be}}%
\thanks{$^{2}$Department of Mathematics, Fine Hall, Princeton University, USA
        {\tt\small Bartv@math.princeton.edu}}%
}

\newcommand{\change}[1]{#1}

\newcommand{\Grad}{\mathrm{Grad}}

\newcommand{\trace}{{\rm Trace}}

\newcommand{\OG}[1]{{\mathcal{O}({#1})}}

\newcommand{\mat}[1]{{\bf #1}}
\newcommand{\subject}{\mathrm{subject\  to}}
\newcommand{\PSD}[2]{\mathrm{S}_+({#1},{#2})}

\begin{document}

\maketitle
\thispagestyle{plain}
\pagestyle{plain}

\begin{abstract}

We propose a Riemannian optimization approach for computing low-rank solutions of the algebraic Riccati equation. The scheme alternates between fixed-rank optimization and rank-one updates. The fixed-rank optimization is on the set of fixed-rank symmetric positive definite matrices which is endowed with a particular Riemannian metric (and geometry) that is tuned to the structure of the cost function. We specifically discuss the implementation of a Riemannian trust-region algorithm that is potentially scalable to large-scale problems. The rank-one update is based on a descent direction that ensures a monotonic decrease of the cost function. Preliminary numerical results on standard small-scale benchmarks show that we obtain solutions to the Riccati equation at lower ranks than the standard approaches.
\end{abstract}


\section{Introduction}
We look at the following low-rank algebraic Riccati equation in $\mat{X} \in \mathbb{R}^{n\times n}$
\begin{equation}\label{eq:lowrank_riccati_equation}
\mat{A}^T\mat{X} + \mat{XA} + \mat{XB}\mat{B}^T {\mat X} = \mat{C}^T \mat{C},
\end{equation}
where $\mat{A} \in \mathbb{R}^{n \times n}$ is a full rank matrix, $\mat{B} \in \mathbb{R}^{n \times p}$, and $\mat{C} \in \mathbb{R}^{s\times n}$ and $s,p\ll n$. The solution of (\ref{eq:lowrank_riccati_equation}) is expected to be symmetric positive semidefinite and low-rank (the rank is $\ll n$) and is of specific interest in a number of disciplines as noted in \cite{simoncini13a, lin13a, benner13a} and references therein. Uniqueness of the solution results from additional assumptions that $-\mat{A}$ and $-\mat{A} - \mat{BB}^T \mat{X}$ are stable, i.e., all their eigenvalues have negative real part \cite{simoncini13a}.

Although conventional solvers for solving (\ref{eq:lowrank_riccati_equation}), e.g., \change{Matlab's} \texttt{ARE} function, exist, they are computationally expensive and cost at least $O(n^3)$. To circumvent the scaling issue in large-scale problems (large $n$), low-rank solvers are of particular interest. Two state-of-the-art approaches that are tuned to the problem (\ref{eq:lowrank_riccati_equation}) are the Kleinman inexact approach \cite[Section~5]{benner13a} and the projection type Krylov subspace approach or the Galerkin Projection approach \cite[Section~3]{benner13a}, \cite{simoncini13a}. The Kleinman inexact approach relaxes the quadratic nature of the equation (\ref{eq:lowrank_riccati_equation}) into a sequence linear Lyapunov equations which are solved approximately by a low-rank Lyapunov solver at every iteration to guarantee a Newton like convergence to the solution. The projection type approach, on the other hand, is based on a series of smaller subspace projections. At each step, the equation (\ref{eq:lowrank_riccati_equation}) is projected onto a smaller subspace which is then solved using standard solver like the \change{Matlab's} \texttt{ARE}. The subspace is incremented in dimension as iterations progress to get accurate solutions. Finally, the rank of the solution is truncated to obtain low-rank solutions. 

In contrast, we view computing a low-rank solution of (\ref{eq:lowrank_riccati_equation}) as computing a solution (critical point) of the problem
\begin{equation}\label{eq:lowrank_riccati_SDP}
\begin{array}{lll}
\min\limits_{\mat {X}\in \mathbb{R}^{n\times n}}  & \frac{1}{4} \| \mat{A}^T\mat{X} + \mat{XA} + \mat{XB}\mat{B}^T {\mat X} - \mat{C}^T \mat{C} \|_F^2 \\
\subject & \mat{X}  \succeq 0,
\end{array}
\end{equation}
where $\mat{X}\succeq 0$ is the positive semidefiniteness constraint and $\| \cdot\|_F$ is the Frobenius norm of a matrix. The quantity $\| \mat{A}^T\mat{X} + \mat{XA} + \mat{XB}\mat{B}^T {\mat X} - \mat{C}^T \mat{C} \|_F$ is referred to as \emph{residual error}. To find a critical point of (\ref{eq:lowrank_riccati_SDP}), we propose a scheme that alternates between fixed-rank optimization and rank-one updates \cite{journee10a} .

The fixed-rank reformulation of (\ref{eq:lowrank_riccati_SDP}) is defined as
\begin{equation}\label{eq:lowrank_riccati_optimization}
\begin{array}{lll}
\min\limits_{\mat{X}\in \mathbb{R}^{n\times n}}  & \frac{1}{4} \| \mat{A}^T\mat{X} + \mat{XA} + \mat{XB}\mat{B}^T {\mat X} - \mat{C}^T \mat{C} \|_F^2 \\
\subject & \mat{X} = \PSD{r}{n},
\end{array}
\end{equation}
where $\PSD{r}{n}$ denotes the set of rank-$r$ symmetric positive semidefinite matrices of size $n\times n$. We tackle the above fixed-rank problem in the framework of Riemannian optimization \cite{edelman98a, absil08a}. The Riemannian optimization framework embeds the constraint into the search space thereby providing an unconstrained optimization setup on the nonlinear search space $\PSD{r}{n}$. A few Riemannian approaches on $\PSD{r}{n}$ are discussed in \cite{journee10a,vandereycken10a,meyer11c} which also list various \emph{ingredients} that enables us to minimize any smooth cost function in a numerically efficient manner. 

A critical component in invoking the Riemannian framework on the search space is the choice of a \emph{Riemannian metric}, a smoothly varying inner product. Selecting the metric leads to equipping the search space with a \emph{Riemannian structure}. \change{Choice of the metric has profound impact on the performance of optimization algorithms.} Tuning the metric then amounts to \emph{preconditioning} the optimization problems by incorporating the Hessian information in the metric. Designing a tuned Riemannian metric on the search space by taking a weighted $L_2$ norm (using a symmetric positive definite part of the Hessian of a specific problem) has been explored in \cite{mishra12a}. This is, for example, done for the low-rank matrix completion problem with good success in \cite{mishra12a, ngo12a}. For the problem (\ref{eq:lowrank_riccati_optimization}) of interest, we exploit the parameterization of the search space discussed in \cite{journee10a} and follow the developments in \cite{mishra12a} to design a novel Riemannian metric in section \ref{sec:Riemannian_metric} for (\ref{eq:lowrank_riccati_optimization}). 

\change{As an algorithm for the fixed-rank problem (\ref{eq:lowrank_riccati_optimization}), we implement a Riemannian trust-region algorithm that has a provably quadratic rate of convergence near the optimum (refer \cite[Chapter~7]{absil08a} for a convergence analysis). In addition, we combine this fixed-rank optimization with a rank increasing outer iteration. The overall algorithm converges superlinearly to a critical point of (\ref{eq:lowrank_riccati_SDP}). This provides a way to keep a tighter control over the rank of the solution while better minimizing the residual error.}

Our main contribution is a meta scheme, shown in Table \ref{tab:algorithm}, for (\ref{eq:lowrank_riccati_equation}) that is based on a novel Riemannian metric (\ref{eq:Riemannian_metric}) for the fixed-rank optimization problem (\ref{eq:lowrank_riccati_optimization}), discussed in section \ref{sec:Riemannian_geometry}. The optimization-related ingredients that are required to implement an off-the-shelf Riemannian trust-region algorithm \cite[Section~7]{absil08a} for (\ref{eq:lowrank_riccati_optimization}) are discussed in section \ref{sec:ingredients_TR}. In particular, the numerical complexity (per iteration) of this approach is discussed in section \ref{sec:numerical_complexity} which shows that our proposed algorithm is potentially scalable to large problems. Preliminary simulations show encouraging results on standard small-scale benchmarks where we obtain lower residual errors at lower ranks with our scheme than the standard approaches.

\section{The proposed Riemannian geometry on the set of fixed-rank symmetric positive semidefinite matrices}\label{sec:Riemannian_geometry}
Any rank-$r$ symmetric positive semidefinite matrix $\mat{X} \in \mathbb{R}^{n\times n}$ is parameterized as 
\[
\mat{X} = \mat{YY}^T,
\]
where $\mat{Y} \in \mathbb{R}_*^{n \times r}$, the set of full column rank matrices of size $n\times r$. This factorization, however, is not unique as $\mat{X}$ remains unchanged under the transformation $\mat{Y} \mapsto \mat{YO}$ for all $\mat{O} \in \OG{r}$, where $\OG{r}$ is the set of matrices of size $r\times r$ such that $\mat{OO}^T = \mat{O}^T\mat{O} = \mat{I}$. The search space is, therefore, identified with the quotient space $\PSD{r}{n} \simeq\mathbb{R}^{n\times r}_* /\OG{r}$. In other words, the search space is the set of equivalence classes $[\mat{Y}] = \{\mat{YO}: \mat{O} \in \OG{r}  \} $, where $\mat{Y} \in \mathbb{R}_*^{n \times r}$.

Based on the above observation, we reformulate the optimization problem (\ref{eq:lowrank_riccati_optimization}) into an optimization problem of a smooth function $\phi:\mathcal{M} \rightarrow \mathbb{R}$ on a quotient manifold \cite{edelman98a, absil08a}
\begin{equation}\label{eq:lowrank_riccati_optimization_manifold}
\begin{array}{lll}
\min\limits_x  & \phi(x) \\
\subject &  x \in \mathcal{M} =  \overline{\mathcal M} /\OG{r},
\end{array}
\end{equation}
where $\overline{\mathcal{M}} = \mathbb{R}_*^{n \times r}$ (the set of full column rank matrices of size $n\times r$) is the computational space, $\OG{r}$ is the set of orthogonal matrices of size $r\times r$. We represent an element of the total space $\overline{\mathcal M}$ by $\bar{x}$ and its corresponding equivalence class by $x$ such that $x = [\bar{x}] $, i.e., $\bar{x}$ has the matrix representation $\mat{Y}$ and $x$ refers to the equivalence class $[\mat{Y}]$. Similarly, the function $\phi:\mathcal{M} \rightarrow \mathbb{R}$ on the quotient manifold $\mathcal{M}$ is induced by the function $\bar{\phi} : \overline{\mathcal M} \rightarrow \mathbb{R}: \bar{x} \mapsto \bar{\phi}({\bar x})$ on the total space $\overline{\mathcal{M}}$. For our case $ \bar{\phi}({\bar x}) = \| \mat{A}^T\mat{YY}^T + \mat{YY}^T\mat{A} + \mat{YY}^T\mat{B}\mat{B}^T \mat{YY}^T - \mat{C}^T \mat{C} \|_F^2 /4$ which is obtained from the cost function in (\ref{eq:lowrank_riccati_optimization}) with the parameterization $\mat{X} = \mat{YY}^T$.

In this section, first, we derive a novel tuned Riemannian metric on $\overline{\mathcal M}$ which induces a Riemannian metric on the quotient space $\mathcal{M}$. Second, the metric structure leads to concrete ideas of implementing a trust-region algorithm on the quotient manifold $\mathcal{M}$ for which we list all the required ingredients. Third, we discuss the per-iteration numerical complexity of the trust-region algorithm and show the potential scalability of the setup.

\subsection{A symmetric positive definite Hessian approximation}
Our tuned metric is based on preconditioning the Hessian by a symmetric positive definite approximation of it. Computation of the symmetric positive definite approximation of the Hessian of the cost function $\bar{\phi}(\mat Y) = \| \mat{A}^T\mat{YY}^T + \mat{YY}^T\mat{A} + \mat{YY}^T\mat{B}\mat{B}^T \mat{YY}^T - \mat{C}^T \mat{C} \|_F^2 /4$ follows from the computation of the gradient $\Grad_{\mat Y} \bar{\phi}$ at $\mat{Y} \in \mathbb{R}_*^{n\times r}$ and its directional derivative in the direction $\mat{Z} \in \mathbb{R}^{n \times r}$. Shown below are few steps that lead to identifying a symmetric positive definite part of the Hessian.
\begin{equation}\label{eq:Hessian_derivation}
\begin{array}{lrll}
 {\rm D} \Grad_{\mat Y} \bar{\phi} [\mat Z]  = (\mat{AA}^T) \mat{Z} (\mat{Y}^T\mat{Y})  \\
\quad + \mat{Z}(\mat{Y}^T\mat{AA}^T\mat{Y}) \\
\quad + (\mat{BB}^T\mat{Y}\mat{Y}^T\mat{Y}\mat{Y}^T \mat{BB}^T) \mat{Z}(\mat{Y}^T\mat{Y}) \\

\quad +  \mat{Z} (\mat{Y}^T\mat{BB}^T\mat{Y} \mat{Y}^T \mat{Y} \mat{Y}^T \mat{BB}^T\mat{Y}) +  {\rm\  other\ terms} \\
\Leftrightarrow {\rm Euclidean\ Hess} [{\rm vec}(\mat Z)] = \underbrace{((\mat{Y}^T\mat{Y}) \otimes \mat{AA}^T)}_{\rm positive\ definite}{\rm vec}(\mat Z)\\
\quad + \underbrace{(\mat{Y}^T\mat{AA}^T\mat{Y} \otimes \mat{I})}_{\rm positive\ definite} {\rm vec}(\mat Z) \\
\quad  + \underbrace{((\mat{Y}^T\mat{Y})\otimes(\mat{BB}^T\mat{Y}\mat{Y}^T\mat{Y}\mat{Y}^T \mat{BB}^T) )}_{\rm positive\ semidefinite} {\rm vec}(\mat Z)  \\

\quad +  \underbrace{((\mat{Y}^T\mat{BB}^T\mat{Y} \mat{Y}^T \mat{Y} \mat{Y}^T \mat{BB}^T\mat{Y})\otimes \mat{I})}_{\rm positive\ semidefinite} {\rm vec}(\mat Z) \\ 

\quad  +  ({\rm other\ negative\ semidefinite\ terms}) {\rm vec}(\mat Z),
\end{array}
\end{equation}
where $\otimes$ denotes the Kronecker product of matrices, ${\rm vec}(\cdot)$ vectorizes a matrix by stacking the columns of the matrix on top of each other, and ${\rm D} \Grad_{\mat Y} \bar{\phi}  [\mat Z]$ is the standard Euclidean directional derivative of $\Grad_{\mat Y} \bar{\phi} $ in the direction $\mat{Z}$, i.e., ${\rm D} \Grad_{\mat Y} \bar{\phi}  [\mat Z] = \lim\limits_{t \rightarrow 0} (\Grad_{\mat Y + t \mat{Z}} \bar{\phi}    - \Grad_{\mat Y} \bar{\phi})/{t} $. It should be stated that when $\mat{B}= 0$, the approximation in (\ref{eq:Hessian_derivation}) extracts the \emph{dominant} component of the Hessian.

\subsection{The proposed Riemannian metric}\label{sec:Riemannian_metric}
Our choice of the metric \change{is} a weighted $L_2$ metric using the above symmetric positive semidefinite approximation of the Hessian. Consequently, the proposed  metric $\bar{g}_{\bar x} : T_{\bar x }\overline{\mathcal M} \times T_{\bar x} \overline{\mathcal M} \rightarrow \mathbb{R}$ on $\overline{\mathcal M}$ at $\bar{x}$ is $\bar{g}_{\bar x} (\bar{\xi}_{\bar x}, \bar{\zeta}_{\bar x}) =  {\rm vec}(\bar{\xi}_{\bar x})  \mathcal{L}{\rm vec}(\bar{\zeta}_{\bar x})$, where $\mathcal{L}$ is the symmetric positive semidefinite approximation of the Hessian derived in (\ref{eq:Hessian_derivation}), ${\rm vec}(\cdot)$ vectorizes a matrix, and $\bar{\xi}_{\bar x},  \bar{\zeta}_{\bar x}$ are any vectors in the tangent space $T_{\bar x} \overline{\mathcal M}$ (the linearization of the manifold $\overline{\mathcal M}$ at $\bar{x}$) at $\bar{x} \in \overline{\mathcal M}$ with the matrix representation ${\bar x} = \mat{Y} \in \mathbb{R}_*^{n\times r}$ and $T_{\bar x}\overline{\mathcal M} = \mathbb{R}^{n\times r}$ \cite[Example~3.6.4]{absil08a}. Equivalently in matrix form, the metric proposed is
\begin{equation}\label{eq:Riemannian_metric}
\begin{array}{lll}
\bar{g}_{\bar x} (\bar{\xi}_{\bar x}, \bar{\zeta}_{\bar x}) = \trace( \bar{\xi}_{\bar x} ^T \mat{A}_1\bar{\zeta}_{\bar x} \mat{M}_1) +  \trace( \bar{\xi}_{\bar x} ^T \bar{\zeta}_{\bar x} \mat{M}_2), \\ 
\end{array}
\end{equation}
where the auxiliary variables $\mat{M}_1 = \mat{Y}^T\mat{Y} \succ 0$, $\mat{A}_1 = (\mat{AA}^T + \mat{BB}^T \mat{YY}^T \mat{YY}^T \mat{BB}^T) \succ 0$, and $\mat{M}_2 =  (\mat{Y}^T \mat{AA}^T \mat{Y} + \mat{Y}^T \mat{BB}^T \mat {Y} \mat{Y}^T \mat{Y} \mat {Y}^T \mat{BB}^T \mat{Y}) \succ 0$ are introduced as shorthand notations. 

Following \cite{journee10a, absil08a}, it can be readily checked that the proposed metric (\ref{eq:Riemannian_metric}) respects the invariance by the group action of $\OG{r}$, induces a Riemannian metric on the quotient space $\mathcal{M}$, and gives a Riemannian submersion structure to $\mathcal{M}$.  Invariance with respect to the action of $\OG{r}$ is critical to define a valid metric on the quotient space $\mathcal{M}$ \cite[Section~3.6.2]{absil08a}. Observe that the metric (\ref{eq:Riemannian_metric}) can be interpreted as a \emph{preconditioner} for the Hessian of $\bar{\phi}$ since the metric is proposed from the Hessian information.

\subsection{Ingredients of a trust-region algorithm on $\PSD{r}{n}$}\label{sec:ingredients_TR}
Once the Riemannian geometry on $\PSD{r}{n} \simeq  \mathcal{M} = \overline{\mathcal M} /\OG{r}$ is decided, it is conceptually straightforward to implement a Riemannian trust-region algorithm on $\PSD{r}{n}$ following \cite{edelman98a,  absil08a}. At each iteration, the trust-region algorithm builds a locally quadratic model around $x$ and minimizes the function in a neighborhood to obtain a candidate search direction \cite[Chapter~7]{absil08a}. In the Riemannian setup, we list the following ingredients that enable us to implement a Riemannian trust-region algorithm \cite[Algorithm~10 and Section~7.5.1]{absil08a}. 
\begin{itemize}
\item Matrix representation of an element $x$ on the quotient manifold $\mathcal{M}$ and its tangent space $T_{x} \mathcal{M}$ at $x$. For $\mathcal{M} = \overline{\mathcal M}/ \OG{r}$, these are identified from the matrix representations in the computational space $\overline{\mathcal M}$. See \cite[Section~3.6.2]{absil08a}.
\item A way to ``move'' on the manifold given a search direction $\xi_x \in T_x \mathcal{M}$. This is accomplished with a retraction mapping \cite[Definition~4.1.1]{absil08a} on $\mathcal{M}$ that maps a tangent vector onto the manifold. The retraction mapping on $\mathcal{M} = \overline{\mathcal M}/ \OG{r}$ is shown in \cite[Example~4.1.5]{absil08a}.
\item Matrix representation of the Riemannian gradient of $\phi$.
\item Matrix representation of a Riemannian connection on the manifold that captures the \emph{covariant derivative} of a vector field $\xi_x$ in the direction of a vector field $\eta_x$. Once the Riemannian connection is defined, the application of Riemannian Hessian along a vector field is directly obtained in terms of the Riemannian connection of the Riemannian gradient along that vector field. 
\end{itemize}

\subsection{Numerical complexity}\label{sec:numerical_complexity}
A carefully study of the ingredients mentioned in section \ref{sec:ingredients_TR} reveals that all operations cost $O(|\mat{A}|r + nr^2 + r^3)$, except the computation of the Riemannian gradient and the Riemannian Hessian operator, \change{where $|\mat{A}|$ is the number of non-zero elements in $\mat{A}$}. Each computation of the Riemannian gradient and the Riemannian Hessian operator needs (only once for every iteration) to solve the linear system
\begin{equation}\label{eq:linear_solve}
\begin{array}{lll}
\mat{A}_1 \xi \mat{M}_1 + \xi \mat{M}_2 = \mat{Z}, \quad \xi \in \mathbb{R}^{ n\times r}
\end{array}
\end{equation}
where $\mat{Z} \in \mathbb{R}^{n\times r}$ is given and the auxiliary variables $\mat{A}_1,\mat{M}_1,\mat{M}_2 \succ 0$ are defined in section \ref{sec:Riemannian_metric}. \change{By means of the generalized eigenvalue decomposition of  $(\mat{M}_1,\mat{M}_2)$, the system (\ref{eq:linear_solve}) can be transformed into $r$ decoupled linear systems.} It has an overall  computational cost of $O({\rm C}_{{\rm solve}}  r + n r^2 + r^3)$, where where ${\rm C}_{\rm solve} $ is the cost of solving shifted systems of $\mat{A}_1$ which has a \emph{sparse + low-rank} structure. In many large-scale problems, ${\rm C}_{\rm solve} $ can be obtained in $O(n)$, although a robust numerical implementation is work in progress. Finally, the computational cost per iteration of the trust-region algorithm is $O({\rm C}_{{\rm solve}}  r + n r^2 + r^3)$.


\section{A meta-scheme for the algebraic low-rank Riccati equation}\label{sec:meta_scheme}
We propose the meta scheme shown in Table (\ref{tab:algorithm}) for (\ref{eq:lowrank_riccati_equation}) that alternates between fixed-rank optimization (with a trust-region algorithm) and rank-one updates \cite{journee10a}. The scheme monotonically decreases the residual error. Convergence of this scheme to a critical point of (\ref{eq:lowrank_riccati_SDP}) is established directly from the analysis in \cite{journee10a}. The rank-one update is based on finding a \emph{descent} search direction that decreases the residual error. If $\mat{Y} \in \mathbb{R}_*^{n\times r}$ is the output of the Riemannian trust-region algorithm that minimizes (\ref{eq:lowrank_riccati_optimization}), then the descent direction is obtained computing the eigenvector $u \in \mathbb{R}^n$ corresponding to the smallest eigenvalue of the \emph{gradient} of the residual $\mat{A}^T\mat{X} + \mat{XA} + \mat{XB}\mat{B}^T {\mat X} - \mat{C}^T \mat{C} $ at $\mat{X} = \mat{YY}^T$, \change{which is of rank $4r +2s$, where $s$ is the number of columns of $\mat{C}$}. \change{The cost of computing the descent direction is $O(|\mat{A}|r + nr^2 + r^3)$, where$|\mat{A}|$ is the cardinality of the matrix $\mat{A}$}. The initialization iterate for the rank $r+1$ subproblem is obtained by the concatenation $[\mat{Y}\ t u] \in \mathbb{R}^{n \times (r +1)}$ where $t >0$ is an appropriate step-size that produces a sufficient decrease in the residual error.  

\begin{table}[t]
\caption{A meta scheme for (\ref{eq:lowrank_riccati_equation}).}
\begin{center}
{\scriptsize
\framebox[3.2in]{
\begin{minipage}[t]{3.0in}
\begin{enumerate}[Step i)]
\item[Given]{
\begin{itemize}
\item Initialize $r$ to $r_0$, say $r_0 =1$.
\item Initialize the iterate $\mat{Y}_0 \in \mathbb{R}_*^{n\times r_0}$.
\end{itemize}
}
\item[Scheme] We alternate between the following two steps until convergence. 
\item[Step i)]{Compute a stationary point $\mat{Y} \in \mathbb{R}_*^{n \times r}$ of the fixed-rank optimization problem (\ref{eq:lowrank_riccati_optimization}) with the Riemannian trust-region algorithm proposed in section \ref{sec:Riemannian_geometry} initialized from $\mat{Y}_0$.}

\item[Step ii)]{Update the rank to $r+1$ and initialize $\mat{Y}_0 = [\mat{Y}  \ tu]$, where $u \in \mathbb{R}^n$ is the descent direction proposed in section \ref{sec:meta_scheme} and $t > 0$ is an appropriate step-size computed by backtracking}
\end{enumerate}
\end{minipage}
}
}
\end{center}
\label{tab:algorithm}
\end{table}

\section{Numerical comparisons}

\begin{figure*}[t]
\center
\subfigure[Effect of choosing a proper metric]{\label{fig:comparisons_precon}
\includegraphics[scale = 0.13]{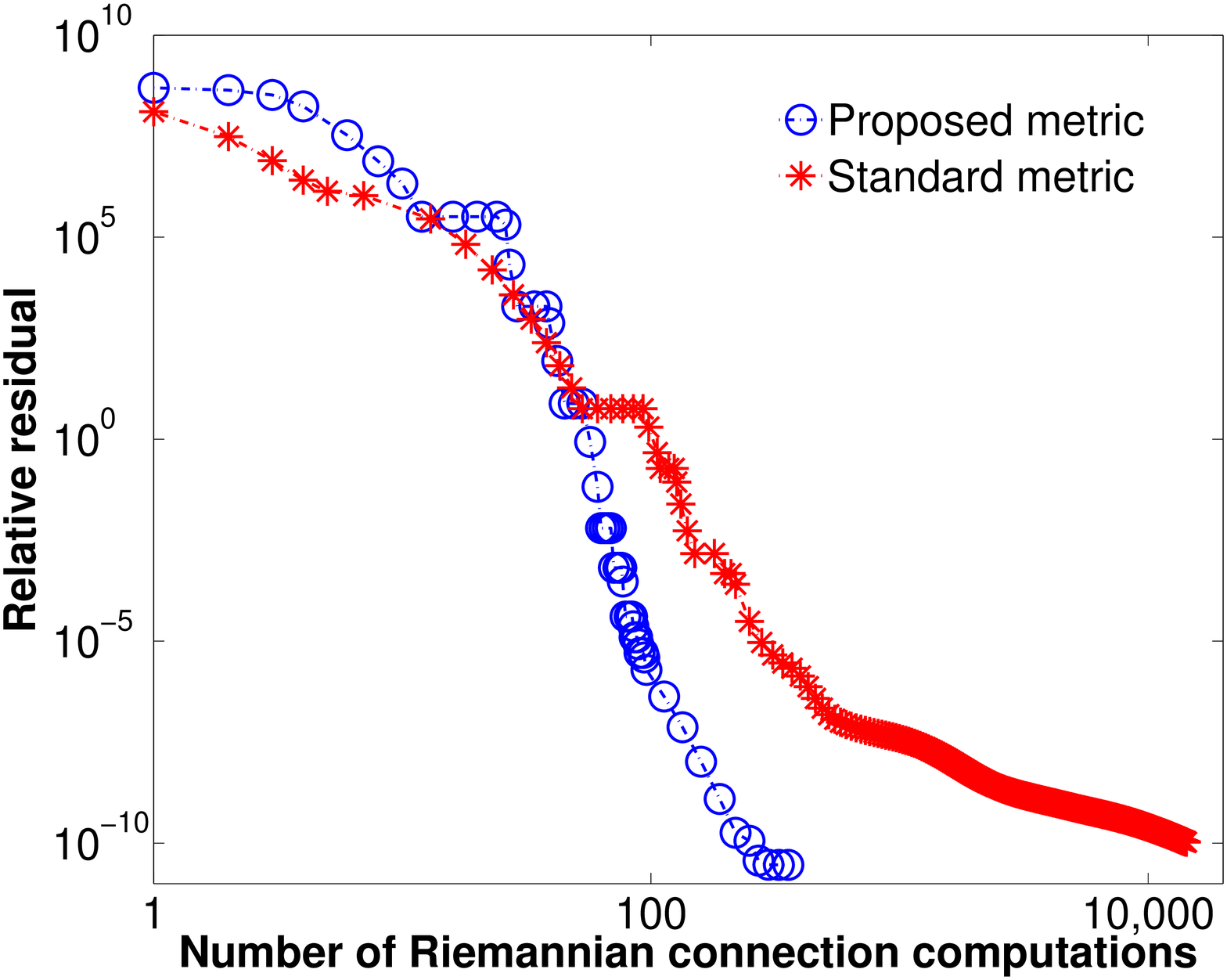}
}
\hspace*{-2em}
\subfigure[Example $1$]{\label{fig:comparisons_laplace}
\includegraphics[scale = 0.13]{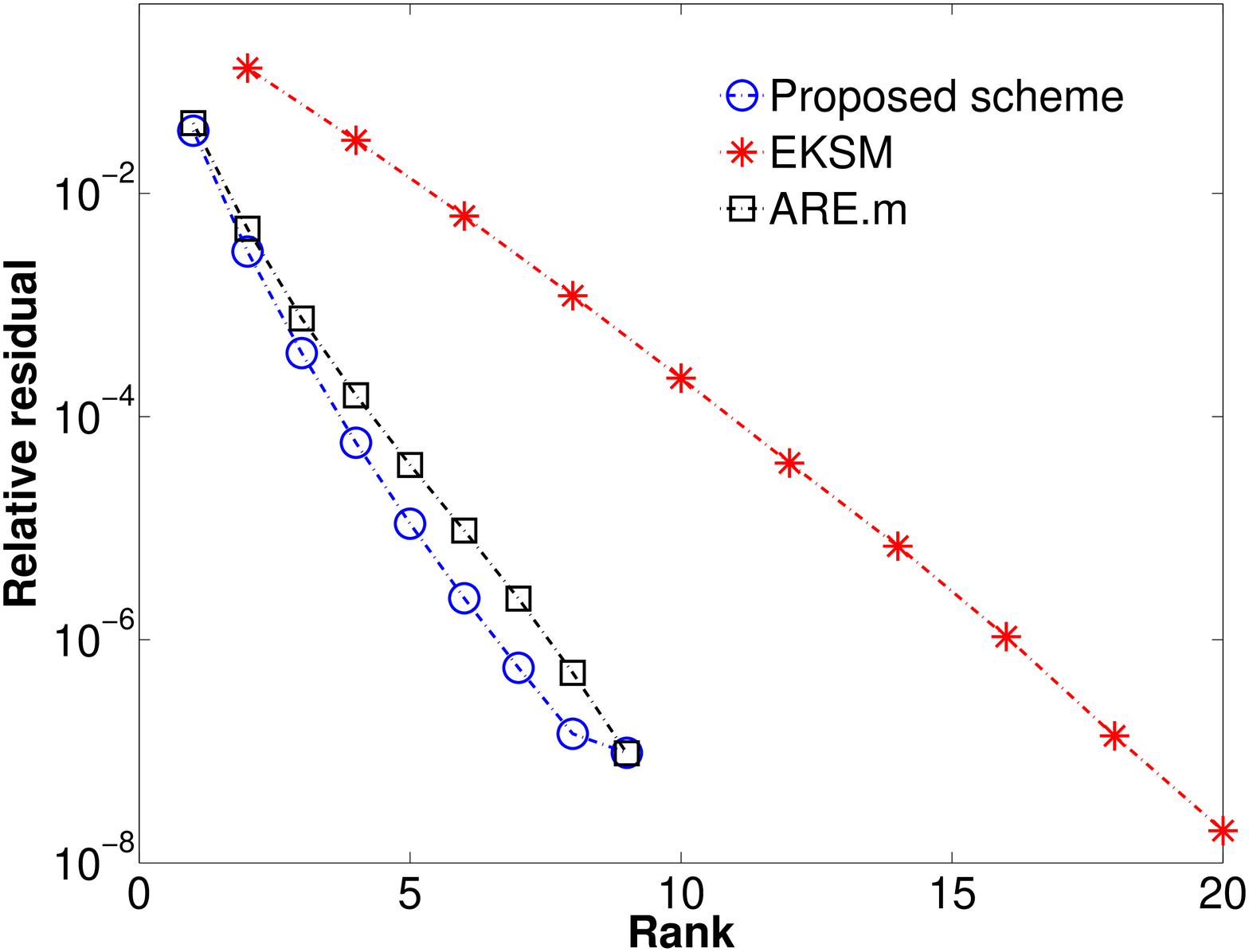}
}
\hspace*{-2em}
\subfigure[Example $2$]{\label{fig:comparisons_toeplitz}
\includegraphics[scale = 0.13]{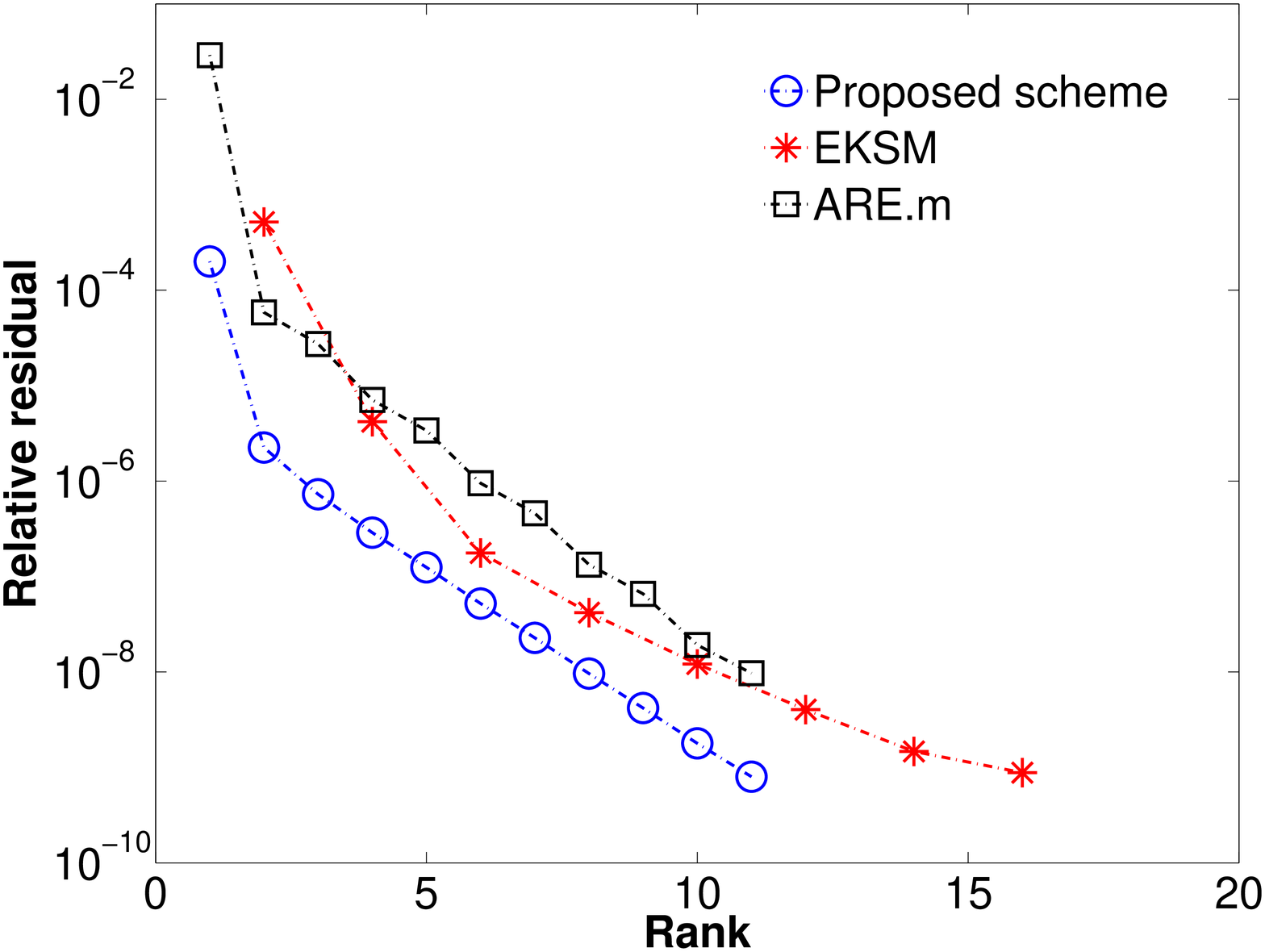}
}
\hspace*{-2em}
\subfigure[Example $3$]{\label{fig:comparisons_heat}
\includegraphics[scale = 0.13]{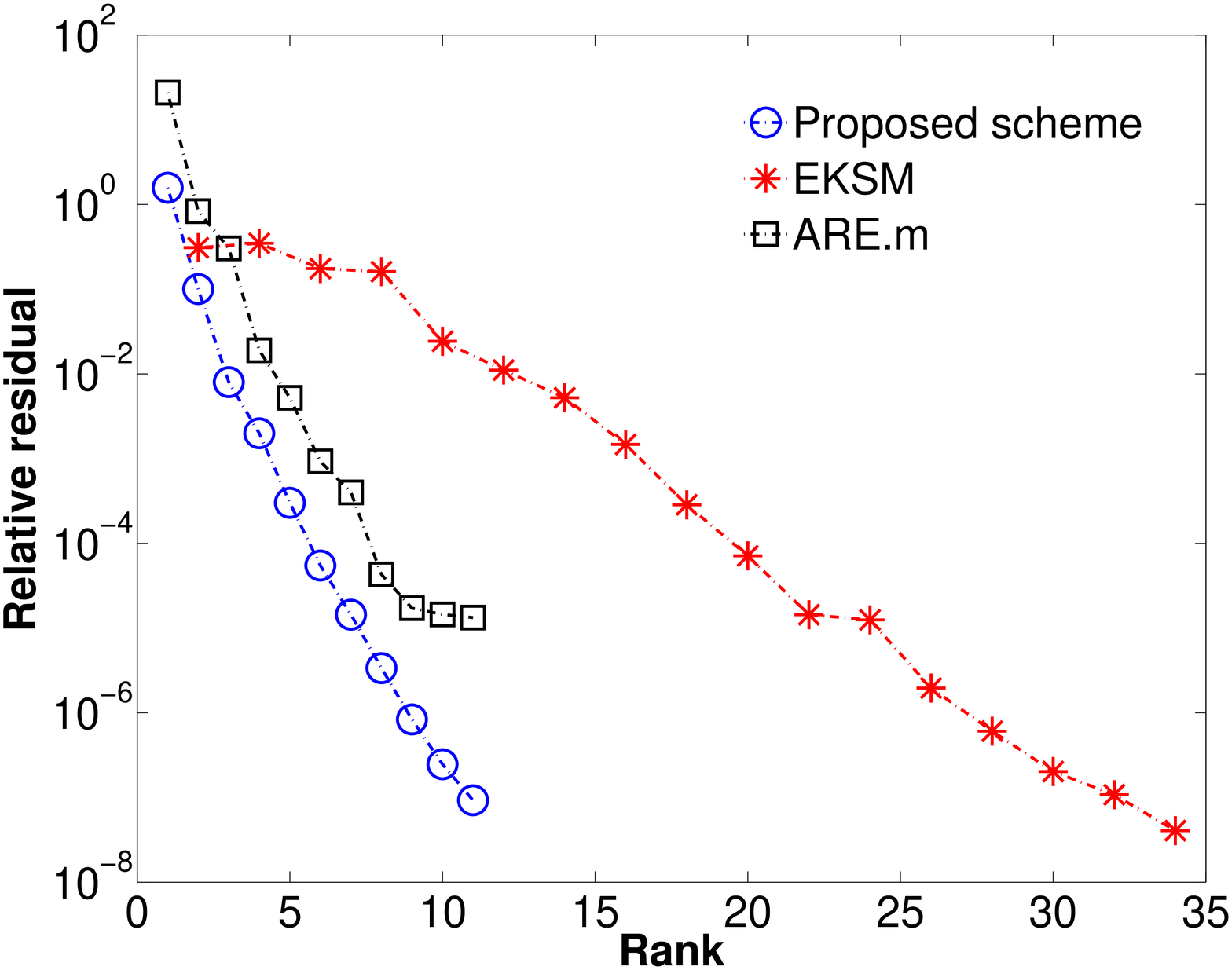}
}
\caption{The benefits of having a tailored Riemannian metric and an optimization-based approach to (\ref{eq:lowrank_riccati_equation}). The proposed scheme in Table \ref{tab:algorithm} leads to smaller relative residual errors at lower ranks for each of the three considered examples. The values on the \texttt{ARE} plot at different ranks are obtained by taking the truncated singular value decomposition of the solution from the Matlab function \texttt{ARE}.}
\label{fig:comparisons}
\end{figure*}
 
To illustrate the notions presented in the paper, we implement our proposed scheme of Table \ref{tab:algorithm} on the open-source Matlab toolbox Manopt \cite{boumal14a}. The toolbox provides an off-the-shelf implementation of the Riemannian trust-region algorithm. Our implementation is available from {\url{http://www.montefiore.ulg.ac.be/~mishra/codes/Riccati.html}}. 
As stopping criteria for our scheme, the fixed-rank optimization is stopped when the norm of the Riemannian gradient norm is below $10^{-10}$ and the rank-one updating is stopped when the relative residual $\| \mat{A}^T\mat{X} + \mat{XA} + \mat{XB}\mat{B}^T {\mat X} - \mat{C}^T \mat{C} \|_F / \|\mat{C}^T \mat{C}\|_F$ is less than $10^{-7}$. The relative residual is computed efficiently by exploiting the low-rank structure \cite[Section~5.3]{lin13a}. For the Riemannian trust-region algorithm, we also limit the number of inner iterations (to solve the trust-region subproblem) to $30$ and the number of outer iterations to $500$.

We first show the effectiveness of the proposed Riemannian metric (\ref{eq:Riemannian_metric}) vis-a-vis the choice of the standard (not tuned) metric in \cite{journee10a} which is $\bar{g}_{\bar x} (\bar{\xi}_{\bar x}, \bar{\zeta}_{\bar x}) = \trace(\bar{\xi}_{\bar x}^T \bar{\zeta}_{\bar x})
$, where $\bar{\zeta}_{\bar x},\bar{\xi}_{\bar x}$ are any vectors in the tangent space at $\bar{x} \in \overline{\mathcal M}$. To this end, we compare the Riemannian trust-region implementations for (\ref{eq:lowrank_riccati_optimization}) on a smaller scale Riccati equation corresponding to \cite[Example~7.1]{lin13a} with $(n, r) = (100,5)$. Figure \ref{fig:comparisons_precon} clearly shows the benefits of a tuned metric that results in a fewer number of Riemannian connection computations (note the log-log scale). The conclusion remains the same across other instances.

We also compare our scheme with the Matlab function \texttt{ARE} and state-of-the-art EKSM algorithm of \cite{simoncini13a} that is based on the Galerkin Projection method. Three examples are considered. Example $1$ corresponds to a smaller scale, $n = 400$, instance of \cite[Example~7.1]{lin13a}. Example $2$ corresponds to the un-normalized case of \cite[Example~7.3]{lin13a} with $n = 500$. Finally, Example $3$ is similar to Example $1$ except $n=500$ and the matrix $\mat{A}$ is a tridiagonal matrix with diagonal entries equal to $2$ and the off-diagonal entries equal to $-1$ that is derived from the discretization of one-dimensional heat equation with Dirichlet boundary. Figures \ref{fig:comparisons_laplace}, \ref{fig:comparisons_toeplitz}, and \ref{fig:comparisons_heat} show that the proposed scheme leads to smaller residual errors at lower ranks for each of the considered examples. With respect to EKSM, it should be stated that the proposed scheme is not competitive in terms of timing. The reason for this is that we traverse through all the ranks one by one minimizing the residual at each rank. This iterative process, while it results in smaller residual errors, is computationally more intensive than EKSM.  

\section{Conclusion}
We have discussed a Riemannian optimization-based approach to low-rank algebraic Riccati equation (\ref{eq:lowrank_riccati_equation}). It leads to the scheme in Table \ref{tab:algorithm} that alternates between fixed-rank optimization and rank-one updates, converging to a critical point of (\ref{eq:lowrank_riccati_SDP}). The fixed-rank optimization problem (\ref{eq:lowrank_riccati_optimization}) is solved with a Riemannian trust-region algorithm on the set of rank-$r$ symmetric positive definite matrices endowed with a novel tuned Riemannian metric (\ref{eq:Riemannian_metric}) that can be seen as a preconditioner for the Riemannian optimization problem. Limited preliminary investigation shows that our approach results in smaller residual errors at lower ranks on standard problem instances. Finding suitable cost functions for (\ref{eq:lowrank_riccati_SDP}) and extending the analogy to other matrix equations will be a topic of future research, as well as having a competitive numerical implementation.

\section*{Acknowledgment}
We thank Valeria Simoncini for providing the codes for the Galerkin subspace methods.


\bibliographystyle{amsalpha}
\bibliography{MV13_Lowrank_Riccati_equation}

\end{document}